\documentclass[11pt]{article}
\usepackage{amsfonts}
\usepackage{latexsym}
\makeatletter

\def\@thm#1#2{\setcounter{#1}{\value{equation}}\refstepcounter{#1}
\stepcounter{equation}
\@ifnextchar[{\@ythm{#1}{#2}}{\@xthm{#1}{#2}}}

\@addtoreset{equation}{section}
\def\qed{\ifhmode\unskip\nobreak\fi\ifmmode\ifinner\else\hskip.5em\fi\fi
 \hbox{\hskip.5em$\Box$\hskip.1em}}

\newenvironment{proof}{\smallskip{\bf Proof.}}{\qed\smallskip}

\newtheorem{theorem}{Theorem}[section]

\newtheorem{definition}{Definition}[section]

\newtheorem{lemma}{Lemma}[section]

\newtheorem{corollary}{Corollary}[section]

\def\@cite#1#2{[{\bf #1}\if@tempswa , #2\fi]}

\def\sideset#1#2#3{%
 \setbox\@ne\hbox{$\displaystyle{\vphantom{#3}}#1{#3}\m@th$}%
 \setbox\tw@\hbox{$\displaystyle{#3}#2\m@th$}%
 \hskip\wd\@ne\hskip-\wd\tw@\mathop{\hskip\wd\tw@\hskip-\wd\@ne
  {\vphantom{#3}}#1{#3}#2}}

\def\up#1{\raise3pt\hbox{\footnotesize#1}}

\long\def\@makecaption#1#2{
 \vskip 10pt
 \setbox\@tempboxa\hbox{#1 #2}
 \ifdim \wd\@tempboxa >\hsize #1 #2\par \else \hbox
to\hsize{\hfil\box\@tempboxa\hfil}
 \fi}

\def\compose{\mathbin{\raise.5ex\hbox{\hskip.08em\tencirc \char'141}}}
\makeatother

\def\bibsameauth{\leavevmode\vrule height .1ex depth 0pt width 
2.3em\relax\,}

\newcommand{\He}{$\mathbf H$}

\newcommand{\Z}{\mathbb Z}

\newcommand{\C}[1]{{\cal #1}}

\newcommand{\ztt}{\hat \otimes}

\begin{document}
\title{Integral Structures on \He{}-type Lie Algebras}
\author{Gordon Crandall \and J\'ozef Dodziuk\thanks{The second author was supported in part by
NSF Grant \#{}DMS-9704743 and the PSC-CUNY Research Award Program}}

\date{\ }
\maketitle

\thispagestyle{empty} 
\bibliographystyle{plain2}

\section{Introduction}
In this paper we prove that every {\He}-type Lie algebra
\cite{kaplan0,kaplan1,kaplan2}
possesses a
basis with respect to which the structure constants are integers. We are
going to call such a basis \emph{an integral basis.} Existence of
an integral basis implies
via the Mal'cev criterion that all simply connected
{\He}-type Lie groups
contain cocompact lattices. Since the Campbell-Hausdorff formula is
very simple for two-step nilpotent Lie groups we can actually avoid
invoking the Mal'cev criterion and exhibit our lattices in an explicit 
way.

The theory of {\He}-type Lie algebras is related very closely to
the theory of Clifford algebras and Clifford modules 
(cf.\ \cite{hus, law-mich}) and we are going to use the classification
of Clifford modules in our construction.

We briefly recall the definition of \He{}-type Lie algebras and establish
notation and conventions for the sequel. Let $\C{U}$ and $\C{V}$ be two
finite-dimensional inner product spaces over $\mathbb R$ of dimensions $m$ and 
$n$
respectively. Let $J: \C{U} \longrightarrow \mbox{End}(\C{V})$, $z
\mapsto J_z$, be a
linear mapping satisfying
\begin{eqnarray}
\label{om1} |J_z(v)| & = & |z|\,|v| \\
\label{om2} J_z\compose J_z (v) & = & - |z|^2\,v
\end{eqnarray}
for all $z \in \C{U}$ and $v \in \C{V}$.
Such a mapping $J$ is called an orthogonal multiplication.  Because of
(\ref{om2}), by the universal property of Clifford algebras
\cite[Proposition 1.1, Chapter 1]{law-mich}, $J$ extends to an algebra homomorphism 
$\phi$ of the Clifford 
algebra $C(\C{U})$ into $\mbox{End}(\C{V})$ so that $\C{V}$ acquires the
structure of a module over $C(\C{U})$. We will often write
$\phi(\alpha)v = \alpha v$ and call $\alpha v$ the Clifford product of 
$\alpha$ and $v$. 

After squaring, polarization of (\ref{om1}) 
 yields
\begin{eqnarray}
\label{om1p} (J_{z_1}(v),J_{z_2}(v) & = & (z_1,z_2)\,|v|^2 \\
\label{om2p} (J_z (v_1),J_z(v_2)) & = & |z|^2 (v_1,v_2)
\end{eqnarray}
holding for all $z,z_1,z_2 \in \C{U}$ and $v,v_1,v_2 \in \C{V}$. In
order not to complicate notation, we use the same notation for norms
and inner products in $\C U$ and $\C V$.
It follows immediately that \begin{equation}
\label{skew} (J_z (v_1),v_2) = -(v_1,J_z(v_2)),
\end{equation}
i.e.\ $J_z(v)$ is a skew-adjoint endomorphism of $\C{V}$, $J_z^*= - J_z$. 
Now the Lie
algebra structure is defined on $\C{N}=\C{U}\oplus\C{V}$ by requiring
that $\C{U}$ be contained in the center and that the bracket of two
elements of $\C{V}$ belong to $\C{U}$ and satisfy $$
(z,[v_1,v_2]) = (J_z(v_1),v_2).$$
In this way, $\C{N}$ becomes a two-step nilpotent Lie algebra which is
referred to as a Heisenberg-type or {\He}-type Lie algebra. 
Orthonormal bases $z_1,\ldots ,z_m$ and
$v_1,\ldots ,v_n$ of $\C{U}$ and $\C{V}$
respectively give rise to the orthonormal basis 
$z_1,\ldots ,z_m,v_1\ldots ,v_n$ of $\C{N}$. The only nonzero structure
constants for $\C{N}$ with respect to this basis occur among numbers 
$A^k_{i,j}$
defined by 
\begin{equation}
\label{struc} [v_i,v_j] = A^k_{i,j}z_k.
\end{equation}

We are now ready to state our main result.
\begin{theorem}\label{main}
For every {\He}-type Lie algebra $\C{N}=\C{U}\oplus\C{V}$
as above there exist orthonormal bases $z_1,\ldots ,z_m$ and
$v_1,\ldots ,v_n$ of $\C{U}$ and $\C{V}$
respectively so that the structure  constants $A^k_{i,j}$
of the Lie algebra
$\C{N}$  with respect to the basis $z_1,\ldots ,z_m,v_1\ldots ,v_n$
are integers and in fact take values $0,1,-1$.
\end{theorem}

The numbers $A^k_{i,j}$ are clearly equal to
$(J_{z_k}(v_i),v_j)=(z_kv_i,v_j)$, i.e.\ depend only on
the Clifford module structure and the inner product of $\C V$.
According to \cite[Proposition 5.16]{law-mich}, every Clifford module,
i.e.\ a finite dimensional module over $C(\C{U})$, admits an inner product
such that the Clifford multiplication by elements of $\C U \subset
C(\C{U})$ is an orthogonal multiplication. We remark that the Clifford
multiplication is orthogonal if and only if elements of $\C U$ act by
skew-adjoint transformations, i.e.\ if and only if (\ref{skew}) holds.

Thus every Clifford module
gives rise to an {\He}-type Lie algebra and Theorem \ref{main} can be
reformulated as a statement about Clifford modules as follows.

\begin{theorem}\label{main1}
Given an inner product space $\C U$ with an orthonormal basis
$z_1,\ldots , z_m$ and a module $V$ over the Clifford algebra
$C(\C{U})$
there exists an inner product $(\cdot ,\cdot)$ on V and an orthonormal 
basis $v_1\ldots,v_n$ of 
$\C{V}=(V,(\cdot ,\cdot))$ so that the Clifford
multiplication $J_z(v)=zv$ by elements of $\C U$ satisfies (\ref{om1}) and
(\ref{om2}) and $(z_iv_p,v_q)$ is equal to $0,1$ or
$-1$ for all $i,p,q$.
\end{theorem}
We will say that the choice of an inner product and a basis for $V$
as above is an \emph{integral structure} and that $V$ with this additional
structure is an \emph{integral Clifford module}. Observe that the
fact that the values of $(e_iv_p,v_q)$ are $0,1,-1$ is equvalent to
the assertion that each of the generators $e_i$ acts on the basis
$v_1,v_2, \ldots ,v_n$ by a permutation and, possibly, some changes of
sign.

From now on we will abandon the notational distinction between an
inner product space $\C{U}$ and the underlying vector space $U$.

\section{Proof of Theorem \ref{main1}}
Since Clifford modules are completely reducible \cite[p.\ 31]{law-mich},
it suffices to prove Theorem \ref{main1} 
for irreducible Clifford modules.
Suppose $\dim \C{U}=k$ and let $e_1,\ldots , e_k$ be an orthonormal
basis of $\C U$. We are going to identify $\C U$ with ${\mathbb R}^k$
and $C(\C U)$ with 
the algebra $C_k$ generated over $\mathbb R$ by $e_1,\ldots , e_k$
subject to relations
\begin{equation}\label{clifford}
e_ie_j = - e_je_i, \qquad\qquad {e_i}^2 = -1.
\end{equation}

According to the classification of irreducible Clifford modules, for
every $k \neq 3 \pmod{4}$ there exists only one isomorphism class of 
irreducible modules
over $C_k$. If $k=3 \pmod{4}$ there are two such classes, 
but the dimensions as vector spaces over $\mathbb R$ of non-isomorphic
modules 
are equal and the {\He}-type groups associated to them
are isomorphic. We will denote an irreducible module over $C_k$ by
$V_k$.

Classification of irreducible Clifford modules proceeds by induction on
$k$, cf.\ \cite[p.\ 33]{law-mich}, and we shall retrace this induction 
proving 
at every stage existence of an integral basis.  It will be
convenient to first classify ${\mathbb Z}_2$-graded Clifford modules.
We briefly recall their definition. Denote by $C_k^0$ the subspace
generated by products of even numbers of generators $e_1,\ldots , e_k$
and by $C_k^1$ the subspace generated by products of an odd number of 
generators. Then $C_k^0$ is a subalgebra, $C^0_k\oplus C^1_k = C_k$, and
$$C_k^i \cdot C_k^j \subset C_k^{i+j}$$ for $i,j \in {\mathbb Z}_2$.
A finite dimensional space $W$ over $\mathbb R$ is a 
${\mathbb Z}_2$-graded Clifford module if $W = W^0 \oplus W^1$ and
$$
C^i_k \cdot W^j \subset W^{i+j}
$$
with  $i,j \in {\mathbb Z}_2$.

We need to define an analog of integral structure on 
${\mathbb Z}_2$-graded Clifford
modules (we use definitions and notation of \cite[Chapter 11, Sections
4,6]{hus} regarding $\Z_2$-graded tensor products of Clifford algebras
and Clifford modules).
\begin{definition}\label{special}

Suppose $W^0$ is a module over $C_k^0$. A choice of  
inner product $(\,\cdot\, ,\, \cdot \, )$ on $W^0$ and an orthonormal basis 
$w_1,\ldots w_m$ with respect to this inner product is called 
\emph{integral} if the basis elements are
permuted with a possible change in sign by Clifford multiplication by
all double products $e_ie_j$ and $$
(ze_k w,ze_k w) = |z|^2 |w|^2$$
for all $w \in W^0$ and $z \in {\mathbb R}^k$. We will call $W^0$ with the
inner product and a basis satisfying the conditions above
\emph{integral}.
\end{definition}
If $W^0$ is the $0$-component of a $\Z_2$-graded Clifford module $W$ 
then multiplication by $e_k$ is an isomorphism of $W^0$ onto $W^1$.
Therefore we can
transfer the inner product on $W^0$ to $W^1$ and define the inner
product on $W^0\oplus W^1 = W$ by requiring the two summands to be
orthogonal. A simple calculation then shows that the Clifford
multiplication by elements of ${\mathbb R}^k$ is an orthogonal
multiplication and the elements of the basis 
$w_1,\ldots ,w_m, e_kw_1,\ldots e_kw_m$ are permuted by Clifford 
multiplication by every $e_i$ with a possible
change in sign. In particular, Theorem \ref{main1} holds for $W$. In the
sequel we are only going to use inner products such that $W^0 \perp W^1$
and the multiplication by $e_k$ maps $W^0$ onto $W^1$ isometrically.

Our proof of Theorem \ref{main1} will be carried out by showing that
for every $k$ there exists an irreducible ${\mathbb Z}_2$-graded
Clifford module $W_k$ with an integral structure. This will imply that 
\emph{every} (ungraded) irreducible Clifford module is integral. It will 
follow
\emph{a posteriori} that all ${\mathbb Z}_2$-graded Clifford modules are
integral. The main fact in the classification of Clifford modules is
that if $W_k$ and $W_8$ are irreducible $\Z_2$-graded modules over $C_k$
and $C_8$ respectively, then $W_k\ztt W_8$ is an irreducible
$\Z_2$-graded module over $C_k\ztt C_8 \simeq C_{k+8}$.
Our proof of Theorem \ref{main1} will consist of exhibiting an integral
irreducible $\Z_2$-graded Clifford module $W_k$ for $k=1,\ldots ,8$ and 
then showing that $W_k\ztt W_8$ is integral if $W_k$ carried an integral
structure. To handle the low dimensional cases we need the following
general lemma.
\begin{lemma}\label{low-dim}
Suppose $V$ is an integral module over the algebra $C_k^0$. Then
$W = C_k\otimes_{C_k^0} V$ is a $\Z_2$-graded module over $C_k$ with
the $\Z_2$-grading given by $W^0 = 1\otimes V$ and $W^1=e_k\otimes V$.
Thus $W^0$ has integral structure transfered from $V$ via
the isomorphism $v \mapsto 1\otimes v$. In addition, if $V$ is
irreducible as a $C_k^0$-module then $W$ is irreducible as a
$\Z_2$-graded module over
$C_k$.
\end{lemma}
\begin{proof} This is obvious 
(see \cite[Chapter 11, Proposition 6.3]{hus}).
\end{proof}

Our next task is to describe integral structures on Clifford modules
over $C_k$ for $k \leq 8$. For $k \leq 7$, we refer to the description
of (ungraded) Clifford modules given in \cite{barbano}. 

\subsection*{${\mathbf k=1}$} 
$\mathbb R \simeq i\mathbb R$ is acting on $\mathbb C$ by complex
multiplication. The real and imaginary parts are $0$- and $1$-components
of $\Z_2$-grading respectively and the standard inner product together
with $1$ as the basis of the $0$-component give the integral structure.

\subsection*{${\mathbf k=2}$}
Consider the space of quaternions $\mathbb H$ equipped with the standard inner
product. Identify ${\mathbb R}^2$ with the span of $i,j$ and let it
act on $\mathbb H$ by quaternion multiplication. This makes $\mathbb H$
into an irreducible module over $C_2$. In addition, the decomposition 
${\mathbb H} = \mbox{span}\,\{1,k\} \oplus \mbox{span}\,\{i,j\}$ is
a  
$\Z_2$-grading and the basis $1,k$ gives an integral structure.

\subsection*{${\mathbf k=3}$}
In this case the space of quaternions becomes a module over the Clifford 
algebra $C_3$ represented in the algebra of endomorphisms of $\mathbb H$
as the
subalgebra generated by quaternion multiplications by $i$, $j$, and $k$.
Clearly, the standard basis $1,i,j,k$ is permuted with possibly a change
in sign by these endomorphisms so that $\mathbb H$ as a module over
$C_3$ carries an integral structure. However, this module does not have a
natural $\Z_2$-graded structure. We regard $\mathbb H$ as a module
over $C^0_3 \subset C_3$ and apply Lemma \ref{low-dim} to create an
integral $\Z_2$-graded irreducible $C_3$-module structure on
$C_3\otimes_{C_3^0}
{\mathbb H}$.

\bigskip

To cover the remaining low-dimensional cases we will use the algebra of
octonions $\mathbb O$ with its standard generators $1,i_1,\ldots i_7$
and multiplication table given in \cite[Page 448]{jacobson} with $c=-1$.

\subsection*{${\mathbf k=4}$}
This is analogous to the case $k=2$. ${\mathbb R}^4$ is identified with
$\mbox{span}\{i_1,i_2,i_3,i_4\}$ which acts on $\mathbb O$ by octonion
multiplication. The resulting Clifford module is $\Z_2$-graded with
the $0$-component equal to $\mbox{span}\{1,i_5,i_6,i_7\}$ and the
$1$-component $\mbox{span}\{i_1,i_2,i_3,i_4\}$. These bases are
orthonormal with respect to the standard inner product and give rise to
an integral structure. This follows from inspection of the
multiplication table.

\subsection*{${\mathbf k=5,6,7}$}
We treat these three cases simultaneously. The Clifford algebra $C_k$
can be represented as an algebra of endomorphisms of $\mathbb O$
generated by transformations of octonion multiplications by $i_j$, $1
\leq j \leq k$. The resulting Clifford module $V_k$ is irreducible but not 
$\Z_2$-graded. As above, we regard it as a module over $C_k^0$ which 
allows us to create an irreducible, $\Z_2$-graded, integral Clifford module
$W_k = C_k\otimes_{C_k^0} {\mathbb O}$ using Lemma \ref{low-dim}.

\subsection*{${\mathbf k=8}$}
This uses the isomorphism, for every $k\geq 2$, 
$\phi : C_{k-1} \longrightarrow C_k^0$
defined as follows \cite[Chapter 11, Section 6]{hus}. Let 
$x=x_0+x_1$ be the decomposition of $x \in C_k$ into its $0$ and $1$
components. Then $\phi(x)=x_0 + e_k  x_1$ where, given standard
generators $e_1,\ldots ,e_k$ of $C_k$, we regard $C_{k-1}$ as the
subalgebra generated by $e_1, \ldots , e_{k-1}$. The octonions are a
module over $C_7$ as above (we relabel $i_j$ as $e_j$) which allows us
to regard $\mathbb O$ as a module
over $C_8^0$ by defining the multiplication $x o$ as
$\phi^{-1}(x)  o$. Since $\phi^{-1}(e_j)= e_j$ and
$\phi^{-1}(e_8 e_j) =e_j$ for $j\leq 7$ and $\mathbb O$ was
integral as a module over $C_7$, the resulting module over $C_8^0$ is
integral. Applying Lemma \ref{low-dim} we obtain an integral,
irreducible, 
$\Z_2$-graded module $W_8 = C_8 \otimes_{C_8^0} {\mathbb O}$. We remark 
that $W_8$ is irreducible as an ungraded Clifford module as well. This
is because $W_8^0$ is irreducible as a module over $C_8^0$ (since
the pair ($W_8^0, C_8^0)$ is isomorphic to $(V_7,C_7)$ and $V_7$ is
irreducible as a module over $C_7$).

\bigskip

We are now ready for the inductive step in the argument. Suppose $W_k$
and $W_l$ are integral, $\Z_2$-graded Clifford modules over $C_k$ and
$C_l$ respectively. Let $e_1,\ldots ,e_k$ and $f_1,\ldots ,f_l$ be the 
standard generators of $C_k$ and $C_l$. Suppose $W_k^0$ and $W_l^0$ are
equipped with inner products and integral bases $v_1,\ldots ,v_m$ and
$w_1,\ldots ,w_n$ respectively. Recall that as a vector space
$W_k\ztt W_l$ is isomorphic to $W_k\otimes W_l$ so that we can equip
$W_k\ztt W_l$ with an inner product by requiring that $$
(x_1\otimes y_1,x_2\otimes y_2) = (x_1,x_2)\cdot(y_1, y_2).$$
In particular, different components $W_k^i\otimes W_l^j$ for $(i,j)\in
\Z_2\times\Z_2$ are orthogonal. Under the isomorphism $C_k\ztt C_l
\simeq C_{k+l}$, cf. \cite[Corollary 4.8]{hus} (note that our notation
differs from the notation in \cite{hus} since we use the symbol
$\ztt$ instead of $\otimes$ to denote the tensor product of ${\mathbb
Z}_2$-graded algebras), the standard generators of 
$C_{k+l}$ correspond to $e_1\otimes 1, \ldots ,e_k\otimes 1, 1\otimes
f_1, \ldots 1\otimes f_l$. To exhibit an integral structure we need a
basis of $(W_k\otimes W_l)^0 = (W_k^0\otimes W_l^0) \oplus (W_k^1\otimes
W_l^1)$ in addition. The basis that we are going to use is
\begin{equation}\label{basis}
\{v_i\otimes w_j, e_kv_i\otimes f_l w_j \mid 1\leq i\leq m, 1\leq
j \leq n\}.
\end{equation}
\begin{lemma}\label{tensor}
The bases described above give rise to an integral structure on $(W_k\ztt
W_l)^0$ and, consequently, on $W_k\ztt W_l$.
\end{lemma}

\begin{proof}
We first verify that the Clifford multiplication by elements of the
space $U$ spanned by the generators of $C_k\ztt C_l$ is orthogonal.
Recall that it suffices to verify that each element $z\in U$ acts on
$W_k\ztt W_l$ as a skew-symmetric endomorphism. We abuse the notation
and write $z$ for an endomorphism associated with $z$. Clearly,
\begin{eqnarray*}
(e_i\otimes 1)^* = e_i^*\otimes 1 = -e_i\otimes 1\\
(1\otimes f_j)^* = 1\otimes f_j^* = -1\otimes f_j
\end{eqnarray*}
since Clifford multiplications on $W_k$ and $W_l$ are orthogonal. Since
every element of $U$ is a linear combination of such products, our
assertion follows.

To verify properties of the integral structure, it suffices to show that
elements of the basis (\ref{basis}) are permuted up to sign by the
double products of generators. There are three cases to consider
$(e_i\otimes 1)\cdot (e_j\otimes 1)=e_ie_j\otimes 1$, 
$(1\otimes f_i)\cdot (1\otimes f_j)=1\otimes f_if_j$, and
$(e_i\otimes 1)\cdot (1\otimes f_j)=(e_i\otimes f_j)$. In the first two
cases the action is as desired since the multiplication by $e_p$'s
permutes $v_1,\ldots ,v_m$ up to sign and the multiplication by $f_q$'s 
acts the same way on $w_1,\ldots ,w_n$. To treat the third case, note
that $e_p=\pm e_ke_pe_k$ and $f_q=\pm f_lf_qf_l$. Thus, up to signs,
$$
(e_p\otimes f_q)\cdot (v_i\otimes w_j) = e_ke_pe_kv_i\otimes
f_lf_qf_lw_j= e_kv_{i'}\otimes f_lw_{j'}
$$
since double products of generators of $C_k$ and $C_l$ permute up to
sign the distinguished bases of $W_k^0$ and $W_l^0$ respectively.
Similarly, $$
(e_p\otimes f_q)\cdot (e_kv_i\otimes f_lw_j) = e_pe_kv_i\otimes
f_qf_lw_j= e_kv_{i'}\otimes f_lw_{j'}.
$$
This proves the lemma.
\end{proof}

\begin{corollary} Suppose that $W_k$ is an irreducible, integral,
$\Z_2$-graded module over $C_k$. Then $W_k\otimes W_8$ is an irreducible,
integral, $\Z_2$-graded module over $C_k\ztt C_8\simeq C_{k+8}$ 
with the integral structure described above.
\end{corollary}
\begin{proof}
All assertions except irreducibility are contained in the lemma above.
By \cite[Chapter 11, 6.5]{hus} $\dim_{\mathbb R}W_k\otimes W_8 = 
\dim_{\mathbb R}W_k \cdot \dim_{\mathbb R}W_8 = 16 \dim_{\mathbb R}W_k$
is equal to the dimension of an irreducible, $\Z_2$-graded Clifford
module over $C_{k+8}$. It follows immediately that $W_k\otimes W_8$ is
irreducible.
\end{proof}
An easy induction using the explicit description of irreducible, 
$\Z_2$-graded Clifford
modules over $C_{k}, 1\leq k \leq 8$ and the lemma above yields
existence of an irreducible, $\Z_2$-graded Clifford
module $W_k$ over $C_{k}$ for every $k$. 

Some additional work is required to
show that every (ungraded) Clifford module has an integral structure.
Let $a_k$ be the dimension of an irreducible Clifford module over $C_k$,
and let $b_k$ denote the dimension of an irreducible, $\Z_2$-graded
module over $C_k$. It follows from 
\cite[Proposition 6.3, Chapter 11]{hus} that $a_k=b_k/2$ and from basic 
periodicity \cite[6.5, Chapter 11]{hus} that 
$a_{k+8}=16a_k$ and $b_{k+8}=16b_k$. All values of $a_k$ and
$b_k$ can now be computed from Table 1, which summarizes some of the
information about Clifford modules of low dimensions.

\begin{table}
\begin{center}
\begin{tabular}{|c|c|c|}
\hline
$k$ & $a_k$ & $b_k$\\
\hline
1 & 2 & 2\\
\hline
2 & 2 & 4\\
\hline
3 & 4 & 8\\
\hline
5 & 8 & 16\\
\hline
5 & 8 & 16\\
\hline
7 & 8 & 16\\
\hline
8 & 16 & 16\\
\hline
\end{tabular}
\end{center}
\caption{Dimensions of irreducible Clifford modules}
\end{table}

Recall now (cf.~\cite[Chapter 6]{hus}) that if $k=1,2,4,8 \pmod 8$
then the dimensions over $\mathbb R$ of an irreducible Clifford module
and an irreducible $\Z_2$-graded Clifford module are equal. 
Thus $\Z_2$-graded
Clifford modules $W_k$ constructed above are irreducible as (ungraded)
Clifford modules. The classification also says that for these values of
$k$ there exists exactly one isomorphism class of Clifford modules over 
$C_k$ which proves Theorem \ref{main1} for $k=1,2,4,8 \pmod 8$.

Next consider the case $k=3 \pmod 4$. According to the classification
there are two non-isomorphic Clifford modules over $C_k$. We are going to
show that their direct sum is isomorphic as an ungraded Clifford module
to $W_k$ and that the integral structure of $W_k$ induces integral
structures on the two modules. The key role in the proof of this fact is
played by the ``volume element'' $\omega = e_1e_2\ldots e_k$. $\omega$
belongs to the center of $C_k$ and satisfies $\omega^2 =1$, cf.
\cite[Proposition 3.3, Chapter 1]{law-mich}. In addition, the
multiplication by $\omega$ is a symmetric operator on every orthogonal
Clifford module. This can be seen as follows. If $k=4l+3$, then
$$
\omega^*=e_k^*e_{k-1}^*\ldots e_1^* = 
(-1)^{4l+3}e_ke_{k-1}\ldots e_1 =\\
-(-1)^{\frac{(4l+4)(4l+3)}{2}}e_1e_2\ldots e_k = \omega,
$$
since multiplications by $e_p$'s are skew-symmetric and $e_p$'s
anti-commute. Now define $\phi_+$ and $\phi_-$ to be multiplications by
central elements $(1+\omega)/2$ and $(1-\omega)/2$ respectively. 
$\phi_\pm$ are self-adjoint and satisfiy $\phi_\pm^2=\phi_\pm$, i.e. they are
orthogonal projections onto their images. Since $(1+\omega)(1-\omega)=0$
and $\phi_+ + \phi_-=1$ the ranges of these projections are
perpendicular. Since $1\pm\omega$ are central, they are in fact (ungraded) 
Clifford submodules $V_+$ and
$V_-$ of $W_k$.
Since $\omega W_k^0 = W_k^1$, $\phi_\pm$ is
injective on $W_k^0$. A count of dimensions shows that $W_k=\phi_+ W_k^0
\oplus \phi_- W_k^0$ is an orthogonal direct decomposition of $W_k$ into two
Clifford submodules. $\omega$ acts on $V_\pm$ by multiplication by
$\pm1$ so that the two modules are \emph{non-isomorphic}. To conclude our
analysis we exhibit  integral bases of the two modules.  We do the
argument for $V_+$ since the argument for $V_-$ is identical after 
$1+\omega$ is replaced by $1-\omega$. Let $w_1, \ldots, w_m$ be an
integral basis of $W_0$, i.e.\ an orthonormal basis which is permuted up
to sign by multiplications by double products $e_ie_j$. Let 
$$
v_p = \frac{1}{\sqrt{2}}(1+\omega)w_p \in V_+
$$
for $p=1,\ldots ,m$. Since $\omega W_k^0 = W_k^1$ is 
perpendicular to $W_k^0$,
Pythagorean theorem insures that this is an orthonormal basis.
For a generator $e_i$ of $C_k$ we calculate 
$$
(e_i v_p,v_q) = (e_iw_p + e_i\omega w_p, w_q + \omega w_q)/2.
$$
Observe that $e_i\omega w_p$ and $w_q$ belong to $W_k^0$ while 
$e_iw_p$ and $\omega w_q$ are in $W_k^1$. Therefore the inner product
above simplifies to $$
((e_i w_p,\omega w_q) + (e_i\omega w_p,w_q))/2 = (e_iw_p, \omega
w_q),
$$
since $\omega$ is central and selfadjoint. Now $\omega w_q$ can be
expressed as a product of $e_i$ times a product of an even number of
generators of $C_k$. It therefore follows from the definition of an 
integral basis that $e_i \omega w_q$ is up to sign equal to $e_i
w_{q'}$ so that, finally $$
(e_i v_p, e_i v_{q'}) = (v_p,v_{q'})=\delta_{p,q'}$$
which proves that the basis $v_1,\ldots, v_m$ is integral.

The remaining two cases of Theorem \ref{main1} are $k=5,6 \pmod 8$
say $k=8l+r$, with $r$ equal to $5$ or $6$.
By the classification, in either case there is only one isomorphism
class of irreducible modules over $C_k$ and the dimension over 
$\mathbb R$ of an irreducible module in this class is equal to
the dimension of the Clifford module $V_+$ over $C_{4l+7}$ constructed 
above. Clearly, $C_k$ can be regarded as the subalgebra of $C_{4l+7}$
generated by the first $k$ generators so that $V_{+}=V_k$ becomes a 
module over
$C_k$. It is irreducible since its dimension is that of an irreducible
module and, trivially, integral. Theorem \ref{main1} is proved.

We conclude this section with a very simple proof, available now, of
the fact that two non-isomorphic Clifford modules over $C_k$, $k =3 \pmod
4$ give rise to isomorphic \He{}-type algebras.  This is very well known,
cf.~\cite{ricci}, but we give the proof for completeness. 
The structure constants 
of the \He{}-type Lie algebra $U\oplus V_+$ associated to the Clifford 
module $V_+$
have been computed above and are equal to
$(e_iv_p, v_q)=(e_iw_p,\omega w_q)$ (here $U=\mbox{span}\{e_1, \ldots
,e_k\}$). In the analogous
calculation for $V_-$ $\omega$ is replaced by $-\omega$ so that the
corresponding structure constants for the Lie algebra $U\oplus V_-$
(with respect to the basis $e_1, \ldots ,e_k, v'_1,\ldots , v'_m$ 
where $v_p' = (1/\sqrt{2}) (1-\omega) w_p$ are negatives of the structure
constants for $U\oplus V_+$ if we choose the obvious correspondence
between bases of $V_+$ and $V_-$. However, if we deploy the set 
$-e_1, \ldots ,-e_k, v'_1,\ldots , v'_1, \ldots ,v_m'$ as the basis
for
$U\oplus V_-$ instead, the structure constants for the two Lie algebras
are equal so that the algebras are isomorphic.

\section{Examples and applications}
In this section we show that every simply connected Lie group of
Heisenberg type contains a cocompact lattice. We are also going to
calculate the isoperimetric dimension of such groups.

We choose and fix an integral basis as in Theorem \ref{main}.
Since \He{}-type Lie algebras are two-step nilpotent (in particular
nilpotent) the exponential mapping is a global diffeomorphism of the 
Lie algebra and the group. We will therefore use it to identify the
two. In addition the Campbell-Hausdorff formula that expresses the group
multiplication in terms of the Lie bracket takes a particularly simple
form 
\begin{equation}\label{product}
X\cdot Y = X+Y +{\scriptstyle\frac{1}{2}} [X,Y]
\end{equation}
where $X,Y\in N = U\oplus V$  as in the Introduction. Let 
$U_{\mathbb Z}$ be the abelian subgroup of $U$ consisting of 
linear combination with integer coefficients of the generators
$e_1, \ldots ,e_k$ of $C(U)\simeq C_k$. Let $V_{\mathbb Z}$ be the
lattice in $V$ generated by elements of an integral basis.
Since $[V_{\mathbb Z},V_{\mathbb Z}] \subset U_{\mathbb Z}$ our
candidate for the lattice $L$ in the group $(N,\,\cdot\,)$ is the subset
$$
{\scriptstyle\frac{1}{2}} U \oplus V. 
$$ 
Verification of this claim is
straightforward after first checking that $-X$ is the inverse of $X$
with respect to the group operation. Thus $L$ is a subgroup of $N$ and 
it is obviously discreet. The product of an arbitrary element $u+v$ of
$N$ with $w \in V$ is equal to $u +(1/2)[v,w] +v+w$ and has $u +(1/2)[v,w]$
as its $U$ component and $v+w$ as the component in $V$. Given $v$ we
can choose $w \in V_{\mathbb Z}$ so that $v+w$ has all coefficients
in the interval $[0,1]$ when expanded with respect to the distinguished 
basis of $V$. We can then act on $u +(1/2)[v,w]$ by elements $u_1$ of
$(1/2)U_{\mathbb Z}$. Since this action is by ordinary translation
($U_{\mathbb Z}$ is contained in the center of $N$) we can easily
achieve that $u+v+(1/2)[v,w] + u_1$ has its coefficients with
respect to the basis $e_1, \ldots ,e_k$ in the interval $[0,1/2]$. Hence
all coefficients in the expansion of 
$$
(u+v)\cdot (u_1+w) = u+u_1+\frac{1}{2}[v,w] + v+w
$$
with respect to the integral basis chosen lie in $[0,1]$, which proves
the co-compactness of $L$. 

To calculate the isoperimetric dimension of $N$ we will need to
calculate the ranks of groups in the lower central series of $L$.
Using our explicit description of the product in $N$ one verifies 
very easily that the bracket operation of the Lie algebra 
\emph{coincides} with the commutator operation in the group.
We therefore have $[L,L] =[V_{\mathbb Z}, V_{\mathbb Z}] = 
U_{\mathbb Z}$. The second equality requires a proof.  
Observe first that
$$
|(z,[v,w])| = |(zv,w)| \leq |z|\cdot |v| \cdot |w|
$$
 if $z\in U$, $v,w \in V$. It follows that 
 $$ 
 |[v,w]| \leq |v|\cdot |w|.
 $$ 
Now let $e_i$ be one of
the generators of $C_k$. Then $e_i$ induces a permutation of the
basis vectors of $V$ with a possible change in sign. Thus given $v_p$,
$e_kv_p = \pm v_q$ for some $q$. It follows that $|(e_i,[v_p,v_q])| =
|(e_1v_p,v_q)| = 1 $ so that $e_i = \pm [v_p,v_q]$. Thus
each of the generators $e_i$ is a commutator which proves the claim.
We see therefore that
\begin{displaymath}
L/[L,L] \simeq 
\left ({\scriptstyle \frac{1}{2}}U_{\mathbb Z}/
U_{\mathbb Z}\right ) \oplus V_{\mathbb Z} 
\simeq 
{\mathbb Z}_2 \oplus V_{\mathbb Z}.
\end{displaymath}
We remark that there is no ambiguity in the interpretation of the
quotient above since $[L,L] \subset U$ and the group multiplication
by $u \in U$ amounts to the the translation by $u$ by (\ref{product}).

A theorem of Bass \cite{bass} gives the estimate of growth of $L$ in
the word metric with respect to the a finite set of generators. In our
case we can take for example the generating set $S$ given by the vectors of 
an integral basis and their negatives. The
lower central series of $L$ reduces to $L_0=L \supset L_1=[L,L]
\supset L_2 = {0}$. The
ranks of successive quotients are therefore equal to $\dim V$ and 
$\dim U$. If $d=\dim V + 2 \dim U$ then by Bass' theorem
 $L$ has polynomial growth
of degree $d$, i.e.~the number $g(R)$ of distinct elements of $L$ in the
metric ball of radius $R$ in $L$ satisfies $$
g(R) \geq c R^d.
$$
where the constant $c$ depends on the choice of the generating set.
A result of Coulhon and Saloff-Coste \cite{Coul-SC} asserts 
that the Cayley graph of a 
group of polynomial
growth of degree $d$ satisfies a $d$-dimensional isoperimetric inequality.
Since every Lie group has bounded geometry we can invoke Kanai's theorem
\cite{kanai}
to conclude that the \He{}-type group under consideration satisfies an
isoperimetric inequality of the same kind.  More precisely we have the
following.
\begin{theorem} Suppose $N=U\oplus V$ is an \He{}-type group equipped
with a left-invariant metric. There exists a positive constant $c$ such
that for every relatively compact subset F of $N$ with smooth boundary
$\partial F$ we have
$$
\frac{A(\partial F)}{V(F)^{1-1/d}} \geq c,
$$
where $A(\partial F)$ and $v(F)$ denote the $n-1$-dimensional volume of
the boundary of $F$ and $V(F)$ stands for the volume of $F$.
\end{theorem}


\end{document}